\documentclass[english, 16pt]{article}
\makeatletter
\providecommand{\LyX}{L\kern-.1667em\lower.25em\hbox{Y}\kern-.125emX\@}
\makeatother

\usepackage{Loc}
\usepackage{graphics}
\begin{document}

\begin{center}
\linespread{1} \Large{\textbf{Spreading of infectious diseases on
complex networks with non-symmetric transmission probabilities}}\\
\vspace{1in} \Large{\textbf{Britta Daudert}}$^{a}$,
\Large{\textbf{Bai--Lian Li}}$^{a,b}$

\vspace{1in}

$^{a}$ Department of Mathematics\\
University of California of Riverside\\
Riverside, CA 92521--0135, USA, britta@math.ucr.edu\\

$^{b}$ Department of Botany \& Plant Sciences\\
University of California of Riverside\\
Riverside, CA 92521--0124, USA, bai-lian.li@ucr.edu\\
\end{center}

\vspace{.5in}

\Large
\begin{abstract}
\noindent We model the spread of a SIS infection on Small World
and random networks using weighted graphs. The entry $w_{ij}$ in
the weight matrix W holds information about the transmission
probability along the edge joining node $v_i$ and node $v_j$. We
use the analogy between the spread of a disease on a network and a
random walk performed on this network to derive a master equation
describing the dynamics of the process. We find conditions under
which an epidemic does not break out and investigate numerically
the effect of a non-symmetric weight distribution of the initially
infected individual on the dynamics of the disease spread.
\end{abstract}

\vspace{.5in}


\section{Introduction}

\vspace{.1in}\noindent The study of epidemiological models has
been a subject of great interest for many years. The aim is to
model the spread of a particular infectious disease, reproducing
the actual dynamics of the disease and designing strategies to
control and possibly eradicate the infection. Several approaches
to tackle this problem have been undertaken. The majority of
epidemic models are based on a compartmental model in which the
individuals are grouped according to their disease status
\cite{Ba}, \cite{AnMa}. The basic models describe the number of
individuals that are susceptible to ($S(t)$), infected with
($I(t)$) and recovered from ($R(t)$) a particular disease at time
$t$. The difference in responses between individuals, the
influence of the topological structure of the system and many
other complex aspects of the progression of the disease are
neglected in this approach. Although the simplicity of the model
means a loss of information and reality, it enables us to get a
first glimpse of the inner workings of the dynamics of the disease
spread and makes calculations of threshold values and equilibria
possible. The assumptions of this model lead to two standard sets
of differential equations that have provided the foundation of the
majority of mathematical epidemiology:

1) The Susceptible-Infectious-Removed (SIR) model:
\begin{eqnarray}
\frac{dS}{dt}&=& \beta N-\lambda S-\rho S\nonumber\\
\frac{dI}{dt}&=&\lambda S-\delta I-\rho I\nonumber\\
\frac{dR}{dt}&=& \delta I-\rho R\nonumber\\
\end{eqnarray}
2) The Susceptible-Infectious-Susceptible (SIS) model:
\begin{eqnarray}
\frac{dS}{dt}&=& \delta I-\lambda S\nonumber\\
\frac{dI}{dt}&=&\lambda S-\delta I,\nonumber\\
\end{eqnarray}

\noindent where N is the population size, $\beta$ is the birth
rate, $\rho$ is the natural death rate, $\delta$ is the recovery
rate and $\lambda$ is the infection rate. The SIR model is a
suitable model for infectious diseases that confer lifelong
immunity, for example measles or whooping cough. The SIS model is
mainly used to model the spread of sexually transmitted diseases,
such as chlamydia or gonorrhoea, where repeated infections are
common. In these models a random mixing assumption is made: each
individual has a small and equal chance of coming into contact
with any other individual. Many modifications to this basic
approach have been made to account for more heterogeneities. One
approach is to further subdivide the population into
subpopulations, with different mixing rates in these groups. This
means that the parameter $\beta$ in the above equations is
replaced by a matrix, describing the transmission of infection
between different groups. Nevertheless, the random mixing
assumption, at least within the subgroups, remains unchanged. In
reality, however, it is usually the case that the number of
contacts of an individual is much smaller than the population size
and random mixing does not occur so that the above model can only
serve as a relatively crude approximation. Models that incorporate
network structure avoid the need to rely on the random mixing
assumption. They do so by assigning each person a fixed amount of
contacts. Networks thus capture the permanence of interactions. A
network (or graph) is comprised of a set of nodes and a set of
edges. An edge is a connection (or bond) that links two nodes. Not
all nodes in a network are connected directly by one edge. Nodes
that are connected by one edge are called neighbors. We introduce
following network quantities:

\begin{itemize}
\item $V=\{v_i\}$, the set of nodes.\item $E=\{e_{ij}\}$, the set
of edges. $e_{ij}$ is the edge running from node $v_i$ to node
$v_j$.\item $k_i$, the degree of node $v_i$, i.e., the number of
neighbors of node $v_i$.\item $<k>$, the average degree of the
nodes.\item $k_{max}$, the maximum degree found in the
network.\item $P(k)$, the degree distribution, i.e., $P(k)$ is the
percentage of nodes in the network that have degree $k$.
\end{itemize}

\noindent In a network used to model the spread of a SIS
infection, the nodes represent individuals that are either
infected by or susceptible to the disease under consideration.
Edges represent interactions between individuals: The disease can
only spread from one node to the next if there is an edge
connecting them. If we assign a weight $\varrho_{ij}$,
representing the probability of node $v_i$ to infect node $v_j$,
to each edge $e_{ij}$, we obtain a weighted graph. The weight
matrix $\mathbf{W}=\{w_{ij}\}$ is defined as follows:

\begin{equation}w_{ij}= \Bigg\{
\begin{array}{cc}
 \varrho_{ij}, & \mbox{if}\quad node \quad v_i\quad is\quad a \quad neighbor\quad of\quad node\quad v_j \\
 0, & \mbox{else}.\nonumber
\end{array}
\end{equation}

\noindent This matrix gives information about the network
topology, the relations between individuals in it and the
characteristic of the disease under consideration. The use of
networks in Mathematical Epidemiology has grown exponentially
since the middle of the 20th century. In recent times, starting
with the work of Pastor-Satorras and Vespignani , there has been a
burst of activity on understanding the effects of the network
topology on the rate and pattern of the disease spread
[\cite{PaVe1}--\cite{KeEa}]. Amongst others the main network types
studied are the following:

\vspace{.2in} \textbf{1) Regular lattices}

\vspace{.2in}\noindent In regular lattices each vertex is
connected to its $k$ nearest neighbors, to form either rings
(one-dimensional) or grids (two-dimensional) \cite{Wa}.

\vspace{.2 in}

\textbf{2) Random graphs}

\vspace{.2 in}\noindent  The term random graph refers to the
disordered nature of the arrangements of links between different
nodes . Erd$\ddot{o}$s and R$\acute{e}$yni (ER), in their first
paper, proposed a model to generate random graphs with N nodes and
K links \cite{ErRe}. Starting with N disconnect nodes, these
random graphs are generated by connecting couples of randomly
selected nodes, prohibiting multiple connections, until the number
of edges is K.

\vspace{.2in}\noindent An alternative model for ER graphs is
created by connecting each pair of nodes with probability $0<p<1$.
This procedure results in graphs having different amount of edges
present but the two models show strong similarities and coincide
in the limit of large N. random graphs have degree distribution
approximately Poisson with parameter $<k>$.

\vspace{.2in}

\textbf{3) Small World networks }

\vspace{.2in}\noindent The study of dynamical processes over real
networks has pointed out the existence of shortcuts, i.e. links
that connect different areas of the graph, thus speeding up the
communication between otherwise distant nodes. This is known as
the Small World property . It is mathematically characterized by a
relatively short average path length that depends at most
logarithmically on the network size. This property is observed in
a variety of real networks including random graphs. To distinguish
between random and Small World networks, the Small World property
is often associated with the presence of clustering. Watts and
Strogatz have proposed to define Small World networks as networks
having both, a short average path length (like random graphs) and
a high clustering coefficient \cite{WaSt}.

\vspace{.2in}

\textbf{4) Scale free networks }

\vspace{.2in}\noindent  Networks that have a power law degree
distribution: $P(k)\sim k^{-\gamma}$, $2<\gamma<3$ are called
scale free \cite{BaAl}. In comparison, in a random graph, $P(k)$
decays faster than exponentially. In scale free networks we often
find a significant amount of nodes with very high degree. In the
context of disease spread on a graph, these nodes are aptly called
super spreaders. In scale free networks the average degree $<k>$
is no longer a relevant variable and one expects fluctuations in
$<k^2>$ to play an important role.

\vspace{.2in}\noindent In their work, Shirley and Rushton
investigated how the speed of the disease spread is influenced by
certain topological characteristics of the graph  \cite{ShRu}.
They found that an epidemic spreads fastest on a scale free
network, followed by random graphs and is slowest on regular
lattices. Small World graphs lay in between random and regular
graphs. In their model, the transmission probability between nodes
is homogeneous throughout the network, i.e., $w_{ij}=w$. However,
when modelling the spread of a disease on a network one should,
along with complex topological features, take into account
heterogeneity in the intensity strength between nodes: not every
individual is susceptible to infection or capable to infect its
neighbors to the same degree. The complexity in the capacity and
intensity of the connections also plays an important role in other
real networks like scientific collaboration networks,
air-transportation networks, internet clusters and other large
infrastructure systems. Newman showed that weighted networks can
in many cases be analyzed  by using a simple mapping from the
weighted graph to an unweighted multi graph and then applying
standard techniques for unweighted networks \cite{Ne}. Simonsen et
al. studied diffusion on unweighted and weighted networks and
discovered that the eigenvalues of the transfer matrix describing
the process can be used to recover large scale topological
features of the system \cite{Sim1}--\cite{Sim2}. Most recently,
Vasquez introduced a type-network representation of a Small World
weighted graph to take into account the population heterogeneity
in a very general approach \cite{Va}. He obtained a recursive
equation for the probability distribution of the outbreak size as
a function of time and demonstrated that the expected outbreak
size and its progression in time are determined by the largest
eigenvalue of a certain matrix (reproductive number matrix) and
the characteristic distance between individuals.

\vspace{.2in}\noindent In this work, we focus on the effect of
heterogeneities in the transmission probabilities on the dynamics
of the disease spread. We use the strong analogy between the
spread of a disease on the network and a random walk performed on
that network to derive a master equation describing the dynamics
of the process. We find conditions under which an epidemic does
not break out and investigate numerically the effect of a
non-symmetric weight distribution on the dynamics of the disease
spread.

\vspace{.2in}\noindent This paper is organized as follows: In
section $2$, the epidemic models (\textbf{SIS}) is described.
Section $3$ contains the derivation of the master equation
describing the dynamic process. Section $4$ is concerned with
identifying conditions that prevent the outbreak of an epidemic.
In section $5$, we discuss our numerical simulations and their
results. We close with an overview of the main results and a
discussion in section $6$.

\section{The model}

\vspace{.1in}\noindent In this paper we will consider the standard
\textbf{SIS} epidemic models without birth and death on a Small
World graph.

\subsection{The epidemic model}
\noindent At each time $t$, the population is divided into two
categories: susceptibles, $S(t)\geq 0$, and infectious, $I(t)\geq
0$. We normalize so that $S(t)+I(t)=1$ for all times $t$.
Susceptible members are virgin territory for the disease, whilst
infectious members are both infected and capable of infecting
others with whom they are in \textbf{direct} contact, i.e., their
neighbors. After being infected for time $\tau_{Inf}$, an
individual returns to the susceptible class. In this paper we take
$\tau_{Inf}=1$.

\subsection{The Network}
\noindent In our networks, the vertices of the graph represent the
individuals of the population under consideration, the edges
describe the contact patterns between individuals. The adjacency
matrix \textbf{A} of this network is defined by

\begin{equation}a_{ij}= \Bigg\{
\begin{array}{cc}
 1, & \mbox{if}\quad node \quad v_i\quad and\quad v_j\quad are\quad joined\quad by\quad an\quad edge \\
 0, & \mbox{else}\nonumber
\end{array}
\end{equation}

\noindent and gives information about neighbor relations in the
network. Assigning a weight $w_{ij}$ to the edge connecting node
$v_i$ and node $v_j$, representing the probabilities for the
disease to spread from node $v_i$ to node $v_j$, will result in a
weighted graph giving information about

\begin{enumerate}
\item The network topology describing the population and the
relations between individuals in it and \item The characteristics
of the disease under consideration.
\end{enumerate}

\noindent This information is encoded in the weight matrix
\textbf{W}, defined in the previous section.

\noindent The weights play the important role of conveying
susceptibility and transmissibility levels of individuals. For
example, consider nodes $v_i$ and $v_j$ with $w_{ij}$ being
relatively large in comparison to $w_{ji}$. This could either mean
that node $v_i's$ capability to transmit the disease is very large
and $v_j's$ very low or that node $v_i$ is not very susceptible to
transmission of the disease from node $v_j$ while $v_j$ can easily
catch the disease from $v_i$.\\

\section{Random Walks and Epidemics}

\noindent In this section we will describe the model implemented
to simulate the spread of an infectious disease throughout a
population and derive a master equation governing the dynamics of
the system. The design of the model was inspired by the work of
Alves et al \cite{AlMaPo}. There is an strong analogy between the
model of an epidemic on a network and a random walk performed on
this network that can be described as follows:
\newline Suppose we want to follow the spread of an SIS epidemic
on a social network. The individuals of the network are
represented as nodes. At each node we place a certain amount of
random walkers. Each random walker represents the possibility of
an infection to happen: If a walker moves along the edge $e_{ij}$
from the infected node $v_i$ to the susceptible node $v_j$, node
$v_j$ has been infected with the disease. If no walker moves to
$v_j$, the node remains susceptible. The probability of the walker
moving from $v_i$ to $v_j$ is given by the entry $w_{ij}$ of the
weight matrix \textbf{W}. There is an artificial component to this
model: The number of walkers placed at each node at time $t_0$
depends on the length of our experiment and the infection time
$\tau_{inf}$. Since each walker represents the possibility of one
infection, we need to have a sufficient amount of walkers at each
node at time $t_0$ so that we do not 'run out' of walkers (the
ability to infect) at some time $t$ before the end of our
experiment. If we want to model the spread of a SIS infection for
$T_{max}$ time steps, we need to place at least $\lceil
\frac{T_{max}}{\tau_{inf}+1}\rceil k_j$ walkers at node $v_j$.
This is to ensure that, when the node is infected, an infection
can happen along each bond eminating from our node
throughout the duration of the experiment.\\
Let us define following quantities:

\begin{itemize}
\item Let $N=T_{max}\sum_jk_j$ be the total number of random
walkers participating. (Here, the sum is taken over all nodes.)
\item Let $\mathbf{W}=\{w_{ij}\}$ be the matrix of transmission
probabilities described earlier: $w_{ij}$ is the probability of
infected individual $v_i$ to infect susceptible individual $v_j$.
\item Let $\mathbf{I}(t)=\{I_i(t)\}$, where

\begin{equation}I_i(t)= \Bigg\{
\begin{array}{cc}
 1 & \mbox{if}\quad v_i\quad is\quad infected \\
 0 & \mbox{else}.\nonumber
\end{array}
\end{equation}

\noindent $\mathbf{I}(t)$ is called the infection matrix. The
number of non-zero elements of $\mathbf{I}(t)$ is the total number
of infected individuals at time $t$.
\end{itemize}

\noindent We model the spread of the infection through the
population via a random walk performed on the network: \newline At
$t=t_0$, we place $T_{max}k_i$ walkers at each node $v_i$. Each
walker represents the possibility of an infection to happen.
Choosing $T_{max}k_i$ walkers (and hence $N=T_{max}\sum_jk_j$)
ensures that no node looses the ability to infect its neighbors
(provided the node is infected) before the maximum number of time
steps $T_{max}$ (i.e., the end of the experiment) is reached. At
each time step $t$, at a generic node $v_i$, one of two things may
happen:

\begin{enumerate}
\item \textbf{If node $v_i$ is infected}: $k_i$ walkers may each
make one move: The walkers are allowed to move between adjacent
vertices. What edge, out of the possible outgoing ones, a walker
chooses to move along is picked at random with probability equal
to the weight assigned to this directed edge. If the walker moves
to a susceptible neighboring node, that node will be infected with
the disease. Walkers from infected neighboring nodes may enter.
After one time step the status of $v_i$ is reset to susceptible
$(I_i(t-1)=1, I_i(t)=0)$. \item \textbf{If node $v_i$ is
susceptible}: No walkers leave $v_i$ but walkers from neighboring
infected nodes may enter. As soon as a walker enters the
susceptible node $v_i$, it is considered infected $(I_i(t-1)=0,
I_i(t)=1)$. At the next time step $k_i$ walkers of node $v_i$ may
move to infect any of the susceptible neighbors of node $v_i$.
\end{enumerate}

\vspace{.2in}\noindent Let us denote by $\eta_i(t)$ the percentage
of walkers at node $v_i$ at time $t$. This is the quantity that
will give us information about the dynamics of the disease spread.
We note following properties of $\eta_i(t)$:

\begin{itemize}
\item If node $v_i$ is susceptible at time $t$,
$\eta_i(t+1)-\eta_i(t)>=0$ and the larger this difference is, the
more susceptible to infection is the node.\item If node $v_i$ is
infected at time $t$, walkers are leaving but may also enter from
neighboring nodes so $\eta_i(t+1)-\eta_i(t)$ does not give any
useful information. \item The change in walker density from time
$t$ to time $t+1$ for the susceptible node $v_i$ satisfies:
$$0\leq|\eta_i(t+1)-\eta_i(t)|\leq\sum_{m\in Neigh_i} \frac{k_m}{\sum_jk_j}.$$
\noindent If
$$0\leq|\eta_i(t+1)-\eta_i(t)|=\sum_{m\in Neigh_i}
\frac{k_m}{\sum_jk_j},$$ all neighboring nodes are infected and
all walkers from those nodes move to $v_i$. \item The change in
walker density from time $t$ to time $t+1$ for the infected node
$v_i$ satisfies:
$$0\leq|\eta_i(t+1)-\eta_i(t)|\leq \frac{k_i}{\sum_jk_j}.$$
\noindent If
$$0\leq|\eta_i(t+1)-\eta_i(t)|= \frac{k_i}{\sum_jk_j},$$ all $k_i$
walkers leave node $v_i$ and no walkers enter. \item Hence, for a
generic node $v_i$, we have
\begin{eqnarray}
0\leq|\eta_i(t+1)-\eta_i(t)|& \leq & \sum_{m\in Neigh_i}
\frac{k_{max}}{\sum_jk_j}+\frac{k_{max}}{\sum_jk_j}\nonumber\\
&=& \frac{1}{\sum_jk_j}[{k_{max}}^2+k_{max}]
\end{eqnarray}
 \item Nodes that are extremely susceptible to
infection at time $t$ have large $|\eta_i(t+1)-\eta_i(t)|$ values
($>>\frac{k_i}{\sum_jk_j}$). Infected nodes that satisfy
$|\eta_i(t+1)-\eta_i(t)|\approx\frac{k_i}{\sum_jk_j}$, are very
infectuous. Hence $|\eta_i(t+1)-\eta_i(t)|$ quantifies the extend
to which node $v_i$ participates in the spread of the disease
under consideration.\item Since the infection time $\tau_{inf}$ is
one and $|\eta_i(t+1)-\eta_i(t)|\leq\sum_{m\in N_i}
\frac{k_m}{\sum_jk_j}+ \frac{k_i}{\sum_jk_j}$, for $T_{max}\geq
3$, we have \footnote{See appendix A and B}:
$$\eta_i(T_{max})\leq\frac{1}{\sum_mk_m}\Big[(T_{max}-\lceil\frac{(T_{max}-3)}{2}\rceil)k_{max}+\lceil\frac{(T_{max}-2)}{2}\rceil {(k_{max})}^2 \Big]$$
\noindent and for $T_{max}\geq 1$:
$$\frac{1}{\sum_mk_m}\Big[(T_{max}-\lceil\frac{T_{max}}{2}\rceil)k_{max}+\lfloor\frac{T_{max}}{2}\rfloor\Big]\leq\eta_i(T_{max})$$
\item From the derivation of these bounds, we can conclude that a
generic node $v_i$ satisfies:
$$\frac{1}{\sum_mk_m}\Big[(t-\lceil\frac{t}{2}\rceil)k_{max}+\lfloor\frac{t}{2}\rfloor\Big]\leq\eta_i(t)\leq\frac{1}{\sum_mk_m}\Big[(t-\lceil\frac{(t-3)}{2}\rceil)k_{max}+\lceil\frac{(t-2)}{2}\rceil {(k_{max})}^2 \Big]$$
\item Nodes with large $\eta_i(T_{max})$ values have been infected
with great intensity but have transmitted the disease with low
frequency.\item Nodes with small $\eta_i(T_{max})$ values have
rarely been infected themselves but if infected have contracted
the disease to many of their neighbors. \item $\sum_i\eta_i(t)=1$
for all $t$.
\end{itemize}

\subsection{The Master Equation}

\noindent We now derive the equation governing the dynamics of the
disease spread. This equation is called the master equation. The
change in the walker density is the difference between the
relative number of walkers entering, $J_i^{-}(t)$, and leaving,
$J_i^{+}(t)$, the same vertex over the time interval
$t\longrightarrow t+1$. Hence,

\begin{equation}
\eta_i(t+1)=\eta_i(t)+J_i^-(t)-J_i^+(t).\label{4a}
\end{equation}

\noindent For the moment, let us presume that node $v_i$ is
infected. Then the edge current on the directed edge from vertex
$v_i$ to a neighboring vertex $v_j$ is given by:

\begin{eqnarray}
C_{ij}(t)&=&\frac{k_i}{\sum_jk_j}\frac{w_{ij}}{\sum_{m\in
Neigh_i}w_{im}}\nonumber \\
&=&\frac{\eta_i(t_0)}{T_{max}}\frac{w_{ij}}{\sum_{m\in
Neigh_i}w_{im}}.\nonumber \\
\end{eqnarray}

\noindent The edge current is the fraction of walkers moving along
this edge according to the weight distribution emanating from that
node. Note that $\frac{w_{ij}}{\sum_{m\in Neigh_i}w_{im}}$ is the
probability of a walker deciding on the edge from vertex $v_i$ to
vertex $v_j$. Since no walkers are leaving if the node is
susceptible, we must have:

\begin{equation}
J_i^-(t)=\sum_{j\in Neigh_i}I_j(t)C_{ji}\label{4b}
\end{equation}

\noindent and

\begin{equation}
J_i^+(t)=\sum_{j\in Neigh_i}I_i(t)C_{ij}.\label{4c}
\end{equation}

\noindent and upon substitution into equation (\ref{4a}), we
obtain

\begin{eqnarray}
\eta_i(t+1)&=&\eta_i(t)+\sum_{j\in Neigh_i}I_j(t)C_{ji}-\sum_{j\in
Neigh_i}I_i(t)C_{ij}\nonumber \\
&=&\eta_i(t)+\sum_{j\in
Neigh_i}I_j(t)\frac{\eta_j(t_0)}{T_{max}}\frac{w_{ji}}{\sum_{m\in
Neigh_j}w_{jm}}-I_i(t)\frac{\eta_i(t_0)}{T_{max}}.\label{4d}
\end{eqnarray}

\noindent With the matrix $\mathbf{T}$ defined as follow:

\begin{equation}T_{ij}= \Bigg\{
\begin{array}{cc}
 \frac{w_{ji}}{\sum_{m\in Neigh_j}w_{jm}} & \mbox{if}\quad v_j\quad is\quad a \quad neighbor\quad of\quad v_i\\
 0 & \mbox{else},\nonumber
\end{array}
\end{equation}

\noindent we can rewrite equation (\ref{4d}) in matrix notation.

\begin{eqnarray}
\mathbf{\eta}(t+1)&=&\mathbf{\eta}(t)+I(t)\cdot \mathbf{T}\frac{\mathbf{\eta}(t_0)}{T_{max}}-\mathbf{I}(t)\cdot\frac{\mathbf{\eta}(t_0)}{T_{max}}\nonumber \\
&=&\mathbf{\eta}(t_0)+\frac{1}{T_{max}}(\sum_{\tau=t_0}^t
\mathbf{I}(\tau))\cdot
\mathbf{T}\mathbf{\eta}(t_0)-(\sum_{\tau=t_0}^t
\mathbf{I}(\tau))\cdot\mathbf{\eta}(t_0) \nonumber \\
&=&\mathbf{\eta}(t_0)+\frac{1}{T_{max}}(\sum_{\tau=t_0}^t
\mathbf{I}(\tau))\cdot
(\mathbf{T}\mathbf{\eta}(t_0)-\mathbf{\eta}(t_0)).\label{4e}
\end{eqnarray}

\noindent Equation (\ref{4e}) is the desired master equation. This
equation governs the dynamics of the disease spread. We can see
that $\mathbf{\eta}(t+1)$ depends on the initial walker
distribution, the sum of the infection matrices from time $t_0$ to
$t$, which tells us which nodes have not been infected yet and the
matrix $\mathbf{T}$, which is called the transfer matrix.

\section{Conditions for non-outbreak}

\noindent If $\mathbf{\eta}(t+1)\approx \mathbf{\eta}(t_0)$  for
all $t$, an outbreak of the epidemic does not take place. From
equation (\ref{4e}) we can see that this is the case if either

\begin{enumerate}
\item  $\sum_{\tau=t_0}^t \mathbf{I}(\tau)\approx \mathbf{0}$ \\
or \item $\mathbf{T}\mathbf{\eta}(t_0)\approx\mathbf{\eta}(t_0)$.
\end{enumerate}

\noindent In this paper we will concentrate on the first
condition. We investigate how the distribution of incoming and
outgoing weights of the initial node effects the dynamics of the
disease spread. We define a weighted difference $D_i$ and a nodal
entropy $S_i$ to quantify the difference between incoming and
outgoing transmission probabilities of a generic node $v_i$.
There are two ways in which condition 1 can be realized:\\
\noindent 1) The epidemic does not break out: either non or
relatively few neighbors of the initially infected node get
infected by it and they, in turn, do not infect a significant
number of their neighbors until, after a short number of time
steps, there are no infected individuals in the population. This
happens if the majority of the transmission probabilities
stay below the threshold value of the network. 
\\
\noindent 2) Infection of susceptible individuals takes place over
a long period of time but infected individuals stay localized
around the initially infected node.

\vspace{.2in}\noindent To investigate when this may happen, we
introduce the following difference measure: For node $v_i$, we
define

$$D_i=\frac{\sum_{j\in Neigh_i}(w_{ij}-w_{ji})}{k_i};$$

\noindent i.e., $D_i$ describes the difference between outgoing
($w_{ij}$) and incoming ($w_{ji}$) transmission probabilities of
node $v_i$. If $D_i<0$, the sum of the outgoing weights is smaller
than the sum of incoming weights and vice versa if $D_i>0$. Nodes
with $D_i<0$ are more likely to be infected often but infect their
neighbors with low intensity. Nodes with $D_i>0$ are more
infectious but less susceptible to infection. We expect some sort
of localization of the disease spread around the initial node
$v_1$ to occur if $D_1\ll 0$ and the average of the weights of the
other nodes are sufficiently small. Let us also define an entropy
describing the diversity of the differences of incoming and
outgoing weights of a node:

\vspace{.1in}\noindent Let

$$P_i(d)=\quad \%\quad  of\quad edges\quad  e_{ij}\quad  of\quad  node\quad  v_i \quad with
\quad w_{ij}-w_{ji}=d$$.

\vspace{.2in}\noindent Then  $P_i(d)$ is a probability
distribution on the differences of incoming and outgoing weights
of node $v_i$ (i.e., $\sum_dP_i(d)=1$ for all $i$) and we utilize
this distribution by defining the following nodal entropy in the
standard way: At node $v_i$, we let
$$S_i^{D}=-\sum_d P_i(d)log(P_i(d)).$$

\noindent $S_i^{D}$ gives information about the disorder of these
differences: The larger the value of  $S_i^{D}$, the large the
variety of differences in incoming and outgoing weights.

\section{Numerical investigation on a Small World Graph}

\subsection{The network}

 We create a Small World graph in the following way:
 $N$ points lying in the square $[0,1]\times[0,1]$ are randomly selected. These points
 represent our population. We choose the following network
 defining quantities:

\begin{itemize}
\item $r$: The short distance radius \item $p_r$: The probability
of short distance bond formation\item $R$: The long distance
radius \item $p_R$: The probability of long distance bond
formation \item $w_r$: The transmission probability along short
distance bonds. \item $w_R$: The transmission probability along
short distance bonds.
\end{itemize}

\noindent Points lying a distance less than $r$ away from each
other are connected with probability $p_r$. Points lying a
distance more than $R$ away from each other are connected with
probability $p_R$. Transmission of the disease along short
distance bonds occurs with probability $p_r$, transmission along
long distance bonds occurs with probability $p_R$.

\subsection{Sensitivity to initial conditions}

\noindent To see when $\sum_{\tau=t_0}^t \mathbf{I}(\tau)\approx
\mathbf{0}$, let us set up our Small World graph as follows:\\
The population size $N$ is 2000. All transmission probabilities
across bonds, regardless of the nature of the bond,  except for
node one, are held fixed at $3\%$ (i.e., $p_r=p_R$). This value is
chosen because our numerical simulations with homogeneous
transmission probabilities showed that an epidemic does rarely
break out for values around $3\%$ We will look at two different
scenarios to evaluate the influence of long distance connections
in the network:

\begin{enumerate}
\item $p_R$ is chosen small enough so that the average long
distance degree of a node lies at $0.09$ (Small world graph).
\item $p_R=0$, no long distance connections are permitted (random
graph).
\end{enumerate}

\noindent $p_r$ is chosen so that the average short distance
degree is $25$. Node one, the initially infected, has degree
$k_1=20$ with $19$ short distance connections and one long
distance connection in case one and 20 connections in case 2. We
define the participation ratio $PR(t)$ at time $t$ as the
percentage of the population that has been infected at least once
during the time interval $[0,t]$ and the diameter $Diam(t)$ as the
maximum distance from the initially infected node where an
infected individual can be found during the time interval $[0,t]$.
We investigate how the diameter and the participation ratio change
with varying $D_1$ and $S_1^D$ values, i.e. how the non-symmetry
of the transmission probabilities of the initial node and the
extend to which it is not symmetric effect the dynamics of the
epidemic. The diameter and participation ratio are measured at the
last time step $T_{max}$. $T_{max}$ is chosen large enough so that
either the epidemic has spread through all the population,
resulting in a maximal diameter and participation ratio, an
equilibrium is reached (participation ratio and diameter approach
some limit) or the epidemic has died out for some $t<T_{max}$. Our
numerical simulations resulted in figure number 1.


\vspace{.2in}\noindent We can see that the plots look very similar
for graphs without long distance connections and Small World
graphs. As expected, the introduction of long distance connections
between
nodes drives up the diameter and participation ratio.\\
\noindent We note that the plots of the diameter and participation
ratio in both cases are almost symmetric about the line $D_1=0$.
For small values of $S_1^D$ we observe a basin around $D_1=0$:
diameter and participation ratio stay low as we hoped. But note
that, as $S_1^D$ increases, the behavior around $D_1\approx 0$
becomes increasingly erratic and unpredictable. As we move away
from $D_1\approx 0$ regions but stay in areas where $S_1^D$ is
relatively low, the diameter increases. We can see that the
infection reaches the outskirts of our domain $[0,1]\times[0,1]$.
participation ratios climb as high as $40\%$. Again, as $S_1^D$
increases we notice the smoothness of the plots breaking and see
disorder taking over. Looking at the participation ratio plots we
notice an oddity: At $D_1\approx 0.4$ and $S_1^D\approx 0.65$ we
see a sharp spike much higher than any surrounding participation
ratio values. This behavior is surprising and needs further
investigation.

\section{Results and Discussion}

\noindent We have used the analogy between the spread of a disease
on a network and a random walk performed on this network to derive
a master equation describing the dynamics of the process. We found
two conditions under which an epidemic does not break out. One of
these conditions strongly depends on the transmission
probabilities of the initial node. This lead to the consideration
of a non-symmetric weight matrix $\mathbf{W}$. The majority of
research concerning epidemic modeling on networks assumes that the
transmission probability matrix $\mathbf{W}$ is symmetric. This
means that node $v_i$ infects neighboring node $v_j$ with the same
probability as node $v_j$ infects node $v_i$. We have focused on a
more realistic setting where we take into account the
heterogeneity in the intensity strength between nodes, i.e., a
non-symmetric weight matrix $\mathbf{W}$. In particular, we focus
on the initially infected individual showing these
heterogeneities, leaving all other nodes with homogeneous
transmission probabilities. We chose these small enough so that in
case of a fully homogeneous network an epidemic would not break
out. We numerically investigated the effect of the non-homogeneous
weight distribution of the initial node. To quantify the
heterogeneity in transmission probabilities, we defined two nodal
quantities:

\begin{enumerate}
\item The difference $D_i=\frac{\sum_{j\in
Neigh_i}(w_{ij}-w_{ji})}{k_i}$, describes the difference between
outgoing ($w_{ij}$) and incoming ($w_{ji}$) transmission
probabilities of node $v_i$. \item The entropy $S_i^{D}=-\sum_d
P_i(d)log(P_i(d))$,
where\\
\noindent $P_i(d)=\quad \%\quad  of\quad edges\quad e_{ij}\quad
of\quad  node\quad  v_i \quad with \quad
w_{ij}-w_{ji}=d$,\\
\noindent describing the disorder in the differences $d$.
\end{enumerate}

\noindent Through our numerical simulations we obtained surface
profiles of diameter and participation ratio for a variety of
$D_1$ and $S_1^D$ values for graphs with long distance connections
and without. Graphs without long distance connections are just
random graphs. Introduction of long distance connections creates
Small World graphs. We obtained the following results:

\begin{itemize}
\item The quasi-symmetry around $D_1\approx 0$ for low $S_1^D$
values shows that regardless of whether incoming weights of node 1
are larger than outgoing or vice versa, the dynamics of the
disease spread is very similar.\item An increase in $S_1^D$ values
leads to a breakdown of symmetry and unpredictable dynamics. This
breakdown of symmetry occurs at $S_1^D$ values about $0.5$.
\end{itemize}

\noindent In summary, we have shown that the introduction of even
very few non-symmetric transmission probabilities along edges of
an otherwise homogeneous network changes the dynamic of the
disease spread. For a larger disorder in the differences between
incoming and outgoing weights of the initially infected, the
dynamics becomes unpredictable. The threshold value of this
disorder measured by the entropy $S_1^D$ lies at about $0.5$. It
is quite apparent that a realistic model of the spread of an
infectious disease throughout a population should take into
consideration the heterogeneities in transmission probabilities
between individuals. We have taken a step towards this direction
by investigating the effect of the introduction of a single node
with heterogeneous transmission probabilities. Further research
opportunities present themselves. One should investigate the
effect of heterogeneities in more than one node. One should try
and find a measure quantifying the extend of anti-symmetry in the
network as a whole and investigate if there exist certain
threshold values determining if an epidemic breaks out or not. We
also plan to analyze the second condition for a non-outbreak found
in this paper. This will involve a thorough investigation of the
properties of the transfer matrix $T$. Another direction of
research would be to run the simulation on scale--free networks to
see if we find the same symmetry about $D_1=0$, if a breakdown of
this symmetry occurs and if so, what the $S_1^D$ threshold value
is.

\vspace{.2in}\noindent Network analysis in the context of epidemic
modelling has helped us create more realistic settings to
investigate the spread of infectious diseases throughout a
population. The study of non-symmetric weight matrices may prove
useful in creating even more realistic models taking into account
the variability of transmission probability between individuals of
the population and could lead to more accurate threshold values
predicting the outbreak of an epidemic, hence enabling us to
develop better methods to prevent or eradicate the disease at
hand.

\begin{center}
\Large{\textbf{Acknowledgements}}
\end{center}

\begin{center}
The authors would like to thank Dr. Michel L. Lapidus for his
input and proof reading of the document and Dr. Erin Pearse for
many fruitful conversations. We would also like to thank Tanya
Vassilevska and Amit Chakraborty for their suggestions and
comments.
\end{center}

\section{Appendix}
\subsection{A} \textbf{Derivation of lower bound}

$$\frac{1}{\sum_mk_m}\Big[(T_{max}-\lceil\frac{T_{max}}{2}\rceil)k_{max}+\lfloor\frac{T_{max}}{2}\rfloor\Big ]\leq\eta_i(T_{max})$$

\vspace{.2in}\noindent First, we make some observations:

\vspace{.2in}\noindent $\tau_{inf}=1$ means that a generic node
$v_i$ can be infected by its neighbor at most every second time
step. Since a node can only loose walkers when it is infected and
node $v_{in}$ is the only infected individual at time $t_0$, we
conclude that $v_{in}$ is the node that potentially has the least
number of walkers at time $T_{max}$. We also need to assume that
node $v_{in}$ is of maximum degree $k_{max}$. The following chain
of events will result in node $v_{in}$ having the smallest number
of walkers at the end of the experiment:

\begin{itemize}
\item \textbf{Time $t_0$:} $v_{in}$ is infected, $v_j$ is
susceptible.
$$\eta_{in}(t_0)=T_{max}\frac{k_{max}}{\sum_mk_m}.$$
\item \textbf{Time $t_1$:} All $k_{in}=k_{max}$ walkers leave
$v_{in}$.
$$\eta_{in}(t_1)=(T_{max}-1)\frac{k_{max}}{\sum_mk_m,}$$ \noindent
$v_{in}$turns susceptible. \item \textbf{Time $t_2$:} To keep
loosing walkers at the largest rate, the node needs to be infected
at every second time step. We assume that only one walker causes
the infection at each time, keeping the amount of walkers entering
at a minimum.\\
$\eta_{in}(t_2)=(T_{max}-1)\frac{k_{max}}{\sum_mk_m}+\frac{1}{\sum_{m}k_{m}}.$
\item When infected, the node looses all $k_{in}$ walkers
in the next time interval.\\
$\eta_{in}(t_3)=(T_{max}-2)\frac{k_{max}}{\sum_mk_m}+\frac{1}{\sum_{m}k_{m}}.$
\item When susceptible, one walker enters, causing infection.\\
$\eta_{in}(t_4)=(T_{max}-2)\frac{k_{max}}{\sum_mk_m}+2\frac{1}{\sum_{m}k_{m}}...$
\item Hence, for $t=T_{max}$, we obtain:
$$\eta_{in}(T_{max})=\frac{1}{\sum_mk_m}\Big[(T_{max}-\lceil\frac{T_{max}}{2}\rceil)k_{max}+\lfloor\frac{T_{max}}{2}\rfloor\Big].$$
\end{itemize}

\subsection{B} \textbf{Derivation of upper bound}

$$\eta_i(T_{max})\leq \frac{1}{\sum_mk_m}\Big[(T_{max}-\lceil\frac{(T_{max}-3)}{2}\rceil)k_{max}+\lceil\frac{(T_{max}-2)}{2}\rceil {(k_{max})}^2 \Big]$$

\vspace{.2in}\noindent Let $v_j$ be a neighbor of the initially
infected individual $v_{in}$. If we assume that $v_{j}$ and all
its neighbors are of maximum degree $k_{max}$ and notice that only
walkers sitting at infected nodes can move, we conclude that the
most walkers can accumulate at a node $v_j$, neighbor to $v_{in}$.
The following chain of events will result in $v_j$ obtaining the
largest number of walkers over the time span of the experiment:

\begin{itemize}
\item \textbf{Time $t_0$:} $v_{in}$ is infected, $v_j$ is
susceptible.
$$\eta_j(t_0)=T_{max}\frac{k_{max}}{\sum_mk_m}.$$
\item \textbf{Time $t_1$:} All walkers leave node $v_{in}$ and
move to node $v_j$. Node $v_j$ gets infected. Node $v_{in}$ turns
susceptible.
$$\eta_j(t_1)=(T_{max}+1)\frac{k_{max}}{\sum_mk_m}.$$
\item \textbf{Time $t_2$:} All $k_{max}$ walkers leave node $v_j$:
$$\eta_j(t_2)=(T_{max})\frac{k_{max}}{\sum_mk_m},$$
\noindent infecting each neighbor. The node returns to the state
of susceptibility .\item \textbf{Time} $t_3$: Assuming that all
walkers of all infected neighboring nodes will move to $v_j$, we
obtain:
\begin{eqnarray}
\eta_j(t_3)&=&(T_{max})\frac{k_{max}}{\sum_mk_m}+\sum_{i\in
N_j}\frac{k_{i}}{\sum_{m}k_{m}}\nonumber \\
&=&
(T_{max})\frac{k_{max}}{\sum_mk_m}+\frac{{(k_{max})}^2}{\sum_{m}k_{m}},\nonumber
\end{eqnarray}
since we assumed all neighbors of $v_j$ to be of maximum degree.
\item We now assume that step 2 and 3 repeat until the end of the
experiment. \item \textbf{Time} $t_4$:
$$(T_{max}-1)\frac{k_{max}}{\sum_mk_m}+\frac{{(k_{max})}^2}{\sum_{m}k_{m}}.$$
\item \textbf{Time} $t_5$:
$$(T_{max}-1)\frac{k_{max}}{\sum_mk_m}+2\frac{{(k_{max})}^2}{\sum_{m}k_{m}}...$$
\item
$$\eta_j(T_{max})=\frac{1}{\sum_mk_m}\Big[(T_{max}-\lceil\frac{(T_{max}-3)}{2}\rceil)k_{max}+\lceil\frac{(T_{max}-2)}{2}\rceil {(k_{max})}^2 \Big].$$
\end{itemize}

\end{document}